\documentstyle{amsppt}
\baselineskip18pt
\magnification=\magstep1
\pagewidth{30pc}
\pageheight{45pc}
\hyphenation{co-deter-min-ant co-deter-min-ants pa-ra-met-rised
pre-print pro-pa-gat-ing pro-pa-gate
fel-low-ship Cox-et-er dis-trib-ut-ive}
\def\leaderfill{\leaders\hbox to 1em{\hss.\hss}\hfill}
\def\A{{\Cal A}}

\def\H{{\Cal H}}
\def\J{{\Cal J}}
\def\L{{\Cal L}}

\def\TL{{\Cal T}\!{\Cal L}}

\def\idest{i.e.,\ }

\def\th{{\theta}}

\def\l{{\lambda}}

\def\s{{\sigma}}

\def\P{{\widetilde P}}

\def\T{{\widetilde T}}
\def\te{\widetilde t}

\def\b0{\text{\bf 0}}

\def\ra{{\ \longrightarrow \ }}

\def\lan{{\langle}}
\def\ran{{\rangle}}

\def\zed{{\Bbb Z}}

\def\boxit#1{\vbox{\hrule\hbox{\vrule \kern3pt
\vbox{\kern3pt\hbox{#1}\kern3pt}\kern3pt\vrule}\hrule}}
\def\rabbit{\vbox{\hbox{\kern0pt
\vbox{\kern0pt{\hbox{---}}\kern3.5pt}}}}

\def\tableau#1{
        \hbox {
                \hskip -10pt plus0pt minus0pt
                \raise\baselineskip\hbox{
                \offinterlineskip
                \hbox{#1}}
                \hskip0.25em
        }
}

\def\tabCol#1{
\hbox{\vtop{\hrule
\halign{\strut\vrule\hskip0.5em##\hskip0.5em\hfill\vrule\cr\lower0pt
\hbox\bgroup$#1$\egroup \cr}
\hrule
} } \hskip -10.5pt plus0pt minus0pt}

\def\CR{
        $\egroup\cr
        \noalign{\hrule}
        \lower0pt\hbox\bgroup$
}



\def\blank#1#2{
\hbox to #1{\hfill \vbox to #2{\vfill}}
}


\def\strut{\vrule height10pt depth5pt width0pt}

\topmatter
\title Fully commutative Kazhdan--Lusztig cells
\endtitle

\author R.M. Green and J. Losonczy \endauthor
\affil Department of Mathematics and Statistics\\ Lancaster
University\\ Lancaster LA1 4YF\\ England\\ {\it  E-mail:}
r.m.green\@lancaster.ac.uk\\
\newline
Department of Mathematics\\ Long Island University\\ Brookville,
NY  11548\\ USA\\ {\it  E-mail:} losonczy\@e-math.ams.org\\
\endaffil

\abstract
We investigate the compatibility of the set of fully commutative
elements of a Coxeter group with the various types of
Kazhdan--Lusztig cells using a canonical basis for a
generalized version of the Temperley--Lieb algebra.

\bigskip

\noindent {\it Cellules pleinement commutatives de Kazhdan--Lusztig}

Nous \'etudions la compatibilit\'e entre l'ensemble des \'el\'ements
pleinement commutatifs d'un groupe de Coxeter et les divers types de
cellules de Kazhdan--Lusztig, en utilisant une base
canonique pour une version g\'en\'eralis\'ee de l'alg\`ebre de
Temperley--Lieb.

\bigskip

\endabstract

\thanks
The first author was supported in part by a NUF--NAL award from the
Nuffield Foundation.
\endthanks

\endtopmatter

\centerline{\bf To appear in Annales de l'Institut Fourier}

\vfill\eject

\head Introduction \endhead

The fully commutative elements, $W_c$, of a Coxeter group $W$ may
be defined, following \cite{{\bf 17}}, as the set of elements $w$
with the property that any reduced expression for $w$ may be
obtained from any other by a sequence of interchanges of adjacent,
commuting Coxeter generators.  These elements arise naturally in
connection with the generalized Temperley--Lieb algebras defined
in the simply laced case by Fan \cite{{\bf 2}} and in general by
Graham \cite{{\bf 8}}.

Kazhdan and Lusztig \cite{{\bf 14}} have defined partitions of a Coxeter
group $W$ arising from each of three equivalence relations,
$\sim_L$, $\sim_R$ and $\sim_{LR}$.  The corresponding equivalence
classes of $W$ are known as left cells, right cells and two-sided
cells, and it follows easily from the definitions that two-sided
cells are unions of left (respectively, right) cells.

In this paper, we are concerned with the compatibility of the set
$W_c$ with the various Kazhdan--Lusztig cells. More precisely, we
wish to know when $W_c$ is a union of (left, right, or two-sided)
cells. Our most general result is Theorem 2.2.3, where
compatibility of $W_c$ with the (left, right, or two-sided)
Kazhdan--Lusztig cells is shown to be related to several other
conditions having to do with the Temperley--Lieb quotient of the
associated Hecke algebra.

It has been known at least since the theses of Fan and Graham
(\cite{{\bf 2}, Proposition 6}, \cite{{\bf 8}}) that $W_c$ is a
union of two-sided cells in type $A$. In Graham's thesis a
counterexample is given which shows that this is not the case in
type $E$.  (See also Example 2.2.5.) Our main task in \S3 is to
prove that there is full compatibility between $W_c$ and the
Kazhdan--Lusztig cells in type $B$.  We have also verified this
property for Coxeter groups of types $F_4$, $H_3$, and $H_4$.

These results are reminiscent of some work of Fan and Stembridge
\cite{{\bf 4}, \S3}, who showed that in types $A$, $D$, $E$ and
affine $A$, the set $W_c$ is a union of
Spaltenstein--Springer--Steinberg cells.

We point out that our methods of proof are combinatorial and based
on our previous work on IC-type (``canonical'') bases for
generalized Temperley--Lieb algebras \cite{{\bf 11}, {\bf 12}}. In
types $A$, $D$ and $E$, the canonical basis defined in \cite{{\bf
11}} is a cellular basis, as defined by Graham in \cite{{\bf 8},
\S4}.

As a consequence of the results in this paper, it becomes possible
to describe the fully commutative cells in type $B$ very
explicitly by using the diagram calculus for canonical basis
elements in type $B$, which was given by the first author in
\cite{{\bf 10}, Theorem 2.2.5}.

\head 1. Kazhdan--Lusztig bases and cells \endhead

\subhead 1.1 Kazhdan--Lusztig bases \endsubhead

We begin by recalling the well-known basic properties of Hecke
algebras arising from Coxeter systems.  These properties all
follow easily from the results of \cite{{\bf 14}}.

Let $X$ be a Coxeter graph, of arbitrary type, and let $W = W(X)$
be the associated Coxeter group with distinguished set of
generating involutions $S=S(X)$.  Denote by $\,<\,$ the
Bruhat--Chevalley ordering on $W$.  Let $\A = \zed[v, v^{-1}]$, let
$\A^- = \zed[v^{-1}]$ and let $q = v^2$.

We denote by $\H=\H(X)$ the Hecke algebra associated with $W$. As
an $\A$-module, the Hecke algebra has a basis consisting of
elements $T_w$, with $w$ ranging over $W$, that satisfy $$T_s T_w
= \cases T_{sw} & \text{ if } \ell(sw) > \ell(w),\cr q T_{sw} +
(q-1) T_w & \text{ if } \ell(sw) < \ell(w),\cr
\endcases$$ where $\ell$ is the length function on the Coxeter group
$W$, $w \in W$, and $s \in S$. It will be convenient to work with
a second $\A$-basis $\{\T_w : w \in W\}$ for $\H$, which we obtain
by defining $\T_w = v^{-\ell(w)}T_w$.

\definition{Definition 1.1.1}
We denote by $a\mapsto \bar{a}$ the $\zed$-linear involution on
the ring $\A$ that exchanges $v$ and $v^{-1}$.  This can be
extended to a $\zed$-linear involution $h\mapsto \bar{h}$ on $\H$,
by defining $\overline{\sum_{w \in W} a_w T_w} = \sum_{w \in W}
\overline{a_w}\, T_{w^{-1}}^{-1}$.
\enddefinition

The Kazhdan--Lusztig basis may be characterized as follows.

\proclaim{Theorem 1.1.2 (Kazhdan--Lusztig)} There is a unique
$\A$-basis $\{C'_w: w \in W\}$ for $\H$ such that
\itemitem{\rm (i)}{$\overline{C'_w} = C'_w$ for all $w\in W;$}
\itemitem{\rm (ii)}{$C'_w = \sum_{x \in W} \P_{x, w} \T_x$,
where each $\P_{x,w}$ lies in $\A^-$ and satisfies $$\P_{x, w} =
\cases 1 & \text{ if } x = w,\cr 0 & \text{ if } x \not\leq w,\cr
0 \mod v^{-1} \A^- & \text{ if } x < w.\cr
\endcases$$}
\endproclaim

\demo{Proof} This is a restatement of \cite{{\bf 14}, (1.1c)}.  The
relationship between $\P_{x, w}$ and the Kazhdan--Lusztig
polynomials $P_{x, w}$ of \cite{{\bf 14}} is $\P_{x, w} = v^{\ell(x) -
\ell(w)}P_{x, w}.$ \qed\enddemo

Theorem 1.1.2 leads to the following reformulation of \cite{{\bf 14},
Definition 1.2}.

\definition{Definition 1.1.3}
Let $x, w \in W$ satisfy $x < w$.  Define $\mu(x, w) \in \zed$ to
be coefficient of $v^{-1}$ in $\P_{x, w}$ (\idest the coefficient
of $v^{\ell(w) - \ell(x) - 1}$ in the Kazhdan--Lusztig polynomial
$P_{x, w}$). If $x < w$ and $\mu(x, w) \ne 0$, then we write $x
\prec w$.
\enddefinition

\subhead 1.2 Kazhdan--Lusztig cells \endsubhead

The structure constants of $\H$ with respect to the
Kazhdan--Lusztig basis of \S1.1 give rise to various natural
partitions of the group $W$ into ``cells''.

Although these structure constants are subtle, the product of two
Kazhdan--Lusztig basis elements may be computed in important
special cases by appealing to the following well-known formula,
which is implicit in \cite{{\bf 14}, \S2.2}.

\proclaim{Proposition 1.2.1 (Kazhdan--Lusztig)} Let $s, w \in W$
with $s \in S$.  Then $$C'_s C'_w = \cases \displaystyle{C'_{sw} +
\sum_{{x \prec w} \atop {sx < x}} \mu(x, w) C'_x} & \text{ if } sw
> w,\cr (v+v^{-1}) C'_w & \text{ otherwise,}
\endcases$$ where $\mu(x, w)$ is as in Definition 1.1.3. $\,\square$
\endproclaim

This formula motivates the following definitions.

\definition{Definition 1.2.2}
Let $x, w \in W$.  We write $x \leq_L w$ if there is a chain $$ x
= x_0, x_1, \ldots, x_r = w ,$$ possibly with $r = 0$, such that
for each $i<r$, $C'_{x_i}$ occurs with nonzero coefficient in the
linear expansion of $C'_s C'_{x_{i+1}}$ for some $s \in S$ such
that $s x_{i+1} > x_{i+1}$.  (By Proposition 1.2.1, this implies
$s x_i < x_i$.)

This transitive preorder yields an equivalence relation $\sim_L$
on $W$ (where $x \sim_L w$ if and only if $x \leq_L w$ and $w
\leq_L x$) whose equivalence classes are called the {\it left
cells} of $W$.  The preorder $\leq_R$ on $W$ is defined by the
condition $x \leq_R w \Leftrightarrow x^{-1} \leq_L w^{-1}$, and
the preorder $\leq_{LR}$ is that generated by $\leq_L$ and
$\leq_R$. These preorders yield equivalence relations $\sim_R$ and
$\sim_{LR}$ on $W$ whose equivalence classes are called {\it right
cells} and {\it two-sided cells}, respectively.
\enddefinition

\remark{Remark 1.2.3} It is well known that the definition of
$\leq_L$ given above agrees with the original definition in
\cite{{\bf 14}}. This follows from Proposition 1.2.1 and part (a)
of the proof of \cite{{\bf 13}, Proposition 7.15}.
\endremark

\remark{Remark 1.2.4} It is immediate from the construction of the
left (respectively right, two-sided) cells that they are partially
ordered via $\leq_L$ (respectively $\leq_R$, $\leq_{LR}$).
\endremark

\head 2. Generalized Temperley--Lieb algebras \endhead

\subhead 2.1 Canonical bases for generalized Temperley--Lieb
algebras
\endsubhead

Let $X$ be a Coxeter graph, of arbitrary type. Let $\J=\J(X)$ be
the two-sided ideal of $\H$ generated by the elements $$ \sum_{w
\in \lan s, s' \ran}T_w, $$ where $(s, s')$ runs over all pairs of
elements of $S$ that correspond to adjacent nodes in the Coxeter
graph. (If the nodes corresponding to $(s, s')$ are connected by a
bond of infinite strength, then we omit the corresponding
relation.)

\definition{Definition 2.1.1}
The generalized Temperley--Lieb algebra, $\TL=\TL(X)$, is the
quotient $\A$-algebra $\H/\J$.  We denote the corresponding
epimorphism of algebras by $\th : \H \ra \TL$.
\enddefinition

The algebra $\TL$ may be of finite or infinite rank, and may be of
finite rank even when it is the quotient of a Hecke algebra of
infinite rank.  Graham \cite{{\bf 8}, Theorem 7.1} classified the
algebras of finite rank into seven infinite families: $A$, $B$,
$D$, $E$, $F$, $H$ and $I$.

\definition{Definition 2.1.2}
A product $w_1w_2\cdots w_n$ of elements $w_i\in W$ is called {\it
reduced} if $\ell(w_1w_2\cdots w_n)=\sum_i\ell(w_i)$.  We reserve
the terminology {\it reduced expression} for reduced products
$w_1w_2\cdots w_n$ in which every $w_i \in S$.

We call an element $w \in W$ {\it fully commutative} if it cannot
be written as a reduced product $x_1 w_{ss'} x_2$, where $x_1, x_2
\in W$ and $w_{ss'}$ is the longest element of some parabolic
subgroup $\lan s, s'\ran$ such that $s,s'\in S$ do not commute.
This definition is equivalent to the one given in the introduction
(see \cite{{\bf 17}, Proposition 1.1}).

We define the {\it content} of $w\in W$ to be the set $c(w)$ of
Coxeter generators $s\in S$ that appear in some (any) reduced
expression for $w$.

Denote by $W_c = W_c(X)$ the set of all elements of $W$ that are
fully commutative.

Let $t_w$ denote the image of the basis element $T_w \in \H$ in
the quotient $\TL$.
\enddefinition

\proclaim{Proposition 2.1.3} {\rm [{\bf 8}, Theorem 6.2]} The set $\{
t_w : w \in W_c \}$ is an $\A$-basis for the algebra $\TL$. $\,
\square$
\endproclaim

We now recall a principal result of \cite{{\bf 11}}, which establishes
a canonical basis for $\TL$.  This basis is a direct analogue of
the Kazhdan--Lusztig basis in \S1, although the precise
relationship between the two is not immediate.

\definition{Definition 2.1.4}
The involution $h\mapsto \bar{h}$ on $\H$ induces a $\zed$-linear
ring involution on $\TL$ \cite{{\bf 11}, Lemma 1.4}. We use the bar
notation to represent this map, as well: $\overline{\sum_{w \in
W_c} a_w t_w} = \sum_{w \in W_c} \overline{a_w}\,
t_{w^{-1}}^{-1}$.

Let $\L$ be the free $\A^-$-submodule of $\TL$ with basis $\{\te_w
: w \in W_c\}$, where $\te_w = v^{-\ell(w)} t_w$, and let $\pi :
\L \ra \L/v^{-1}\L$ be the canonical projection.  For each $w\in
W$, we denote by $\L_w$ the free $\A^-$-submodule of $\TL$ with
basis $\{\te_x: x\in W_c,\,x\leq w\}$.
\enddefinition

\proclaim{Proposition 2.1.5} {\rm [{\bf 11}, Theorem 2.3]} There exists
a unique basis $\{ c_w : w \in W_c\}$ for $\L$ such that
$\overline{c_w} = c_w$ and $\pi(c_w) = \pi(\te_w)$ for all $w\in
W_c$. $\, \square$
\endproclaim

We call $\{c_w : w \in W_c\}$ the {\it canonical basis} (or the
{\it IC basis}) of $\TL$.  It depends on the $t$-basis, the
involution on $\TL$ from above, and the $\A^-$-lattice $\L$.

\definition{Definition 2.1.6}
If $s \in S$, we write $b_s \in \TL$ for the element $v^{-1} t_s +
v^{-1}$.  The elements $b_s=\theta(C'_s)$ generate $\TL$ as an
$\A$-algebra.

For each $w\in W_c$, it makes sense to define
$b_w=b_{s_1}b_{s_2}\cdots b_{s_n}$, where $s_1s_2\cdots s_n$ is
any reduced expression for $w$.  It is known that the set
$\{b_w:w\in W_c\}$ is an $\A$-basis for $\TL$; we call it the {\it
monomial basis}. (In types $A$, $D$ and $E$, this agrees with
Graham's cellular basis for $\TL$ as defined in \cite{{\bf 8}}.)

For each $w\in W$, we denote by $\L'_w$ the free $\A^-$-submodule
of $\TL$ with basis $\{b_x: x\in W_c,\, x\leq w\}$. We shall study
the $\A^-$-lattices $\L'_w$ in \S3.3.
\enddefinition

\subhead 2.2. Some general results \endsubhead

One of the main obstructions to understanding the relationship
between the Kazhdan--Lusztig basis of $\H(X)$ and the canonical
basis of $\TL(X)$ is that the set $W_c(X)$ may not be compatible
with the two-sided cells. When a particular type of compatibility
is present, the relationship between the two bases becomes
transparent \cite{{\bf 12}, Proposition 1.2.3}. It will be shown
in \S3 that this is the case when $X$ is of type $B_n$.

Before restricting ourselves to type $B_n$, we shall say something
more about the general problem. The following result is helpful in
this context.

\proclaim{Lemma 2.2.1} The set $\{\th(C'_w) : w \in W_c\}$ is an
$\A$-basis for $\TL$.
\endproclaim

\demo{Proof} This follows easily from \cite{{\bf 11}, Lemma
1.5}.\qed\enddemo

\proclaim{Lemma 2.2.2} Let $c \in \L$ be such that $\pi(c) = 0$
and $\overline{c} = c$. Then $c = 0$.
\endproclaim

\demo{Proof} Express $c = \sum_{w \in W_c} a_w c_w$ as an
$\A^-$-linear combination of the canonical basis.  Since $\pi(c) =
0$, we must have $a_w \in v^{-1}\A^-$ for all $w$.  But $$
c=\overline{c} = \sum_{w \in W_c} \overline{a_w}\,c_w, $$ which
implies that all $a_w = 0$. \qed\enddemo

There is compatibility between the set $W_c$ and the
Kazhdan--Lusztig cells when the equivalent conditions of the
following theorem are satisfied.

\proclaim{Theorem 2.2.3} Let $X$ be an arbitrary Coxeter graph,
and maintain the usual notation, e.g., $\J=\J(X)$,
$W=W(X)$, etc. Then the following are equivalent:
\itemitem{\rm(i)}{The ideal $\J$ is spanned by those elements
$C'_w$ that it contains.}
\itemitem{\rm(ii)}{The ideal $\J$ is spanned by the set
$\{C'_w : w \in W \backslash W_c\}$.}
\itemitem{\rm(iii)}{For each $w \in W\backslash W_c$, one has
$\th(C'_w) = 0$.}
\itemitem{\rm(iv)}{If $w \in W$, then $\th(C'_w)\in \L$ and
$$\pi(\th(C'_w)) = \cases
\pi(c_w) & \text{ if } w \in W_c;\cr 0 & \text{ otherwise.}\cr
\endcases$$}
\itemitem{\rm(v)}{For each $w \in W\backslash W_c$, one has
$\th(\T_w)\in v^{-1}\L.$}
\itemitem{\rm(vi)}{The set $W \backslash W_c$ is closed under
$\leq_{L}$ and so is a union of left cells.}
\itemitem{\rm(vii)}{The set $W \backslash W_c$ is closed under
$\leq_{LR}$ and so is a union of two-sided cells.}
\itemitem{\rm(viii)}{The set $W_c$ is closed under $\geq_{LR}$
and so is a union of two-sided cells.}
\endproclaim

\demo{Note} If $\leq_\Lambda$ is a preorder on a set $\Lambda$ and
$\Lambda' \subseteq \Lambda$, then the statement ``$\Lambda'$ is
closed under $\leq_\Lambda$'' means that whenever $\l_1 \in
\Lambda'$ and $\l_2 \in \Lambda$ are such that $\l_2 \leq_\Lambda
\l_1$, we have $\l_2 \in \Lambda'$.
\enddemo

\demo{Proof} We begin with the equivalence of (i), (ii) and (iii).
The statement (i) implies that the set $ \{ \th(C'_w) : \th(C'_w)
\ne 0 \} $ forms an $\A$-basis for $\TL$.  It follows that this
set must equal the set in Lemma 2.2.1, which implies (ii).  It is
clear that (ii) implies (iii).  If (iii) holds then $\J$ must
contain all elements $C'_w$ with $w \in W \backslash W_c$.
However, $\J$ cannot be any bigger than the span of these elements
by Lemma 2.2.1, so (i) follows.

\medskip

\noindent (iii) $\Leftrightarrow$ (iv). Assume (iii) holds. It is
clear that if $w\not\in W_c$, then $\th(C'_w)\in \L$ and
$\pi(\th(C'_w)) = 0$. The remaining case is dealt with by the
proof of \cite{{\bf 12}, Proposition 1.2.3}. Now assume (iv) and
let $w \not\in W_c$. Since the map $\th$ is compatible (by
\cite{{\bf 11}, Lemma 1.4}) with the involutions on $\H$ and $\TL$
from above, we see that $\overline{\th(C'_w)}=\theta({C'_w})$.
Since $\pi(\th(C'_w)) = 0$, Lemma 2.2.2 gives (iii).

\medskip

\noindent (iv) $\Leftrightarrow$ (v). Note that if $w \in W_c$,
then $\pi(\th(\T_w)) = \pi(\te_w) = \pi(c_w)$ by the definition of
the canonical basis. The equivalence of (iv) and (v) now follows
from the fact that, relative to some total refinement of the
Bruhat--Chevalley order, the (possibly infinite) change of basis
matrices between the basis $\{C'_w : w \in W\}$ and the basis
$\{\T_w : w \in W\}$ are upper triangular with ones on the
diagonal, and all the entries above the diagonal lie in $v^{-1}
\A^-$.

\medskip

\noindent (ii) $\Rightarrow$ (vi). By Definition 1.2.2, it is
enough to check that if $w \not\in W_c$ and $s \in S$ with $sw >
w$, then all the terms occurring in the expansion of $C'_s C'_w$
in Proposition 1.2.1 are parametrized by elements $x \not \in
W_c$. By (ii), $w \not\in W_c$ implies $C'_w \in \J$. Since $\J$
is an ideal, $C'_s C'_w \in \J$. Another application of (ii)
completes the proof.

\medskip

\noindent (vi) $\Rightarrow$ (vii). Suppose $w \not\in W_c$ and $x
\leq_R w$, meaning that $x^{-1} \leq_L w^{-1}$. By symmetry of the
definition of $W_c$, we have $w^{-1} \not\in W_c$, and (vi) shows
that $x^{-1} \not\in W_c$, meaning that $x \not\in W_c$.  It
follows that $W \backslash W_c$ is closed under $\leq_R$. Since $W
\backslash W_c$ is closed under $\leq_L$ and $\leq_R$, it is
closed under $\leq_{LR}$ and is therefore a union of two-sided
cells.

\medskip

\noindent (vii) $\Rightarrow$ (ii). Since $W\backslash W_c$ is
closed under $\leq_{LR}$ and $\{C'_s : s \in S\}$ is a set of
algebra generators for $\H$, it follows that $\{C'_w : w \not\in
W_c\}$ spans a two-sided ideal of $\H$. The generators of $\J$
(see \S2.1) are of the form $v^{\ell(w)}C'_w$ for certain $w
\not\in W_c$, so $\J$ is contained in this ideal. On the other
hand, $\J$ contains $\{C'_w : w \not\in W_c\}$ by Lemma 2.2.1.
Thus, condition (ii) holds.

The equivalence of (vii) and (viii) is obvious. \qed\enddemo

\example{Example 2.2.4}  Consider the finite dihedral case, or in
other words, take $X=I_2(m)$ for $m<\infty$.  Let $w_0$ denote the
longest element of $W(I_2(m))$.  The ideal $\J(I_2(m))$ is spanned
by the single Kazhdan--Lusztig basis element $C'_{w_0}$.  Thus,
the equivalent conditions of Theorem 2.2.3 hold for finite dihedral
groups.
\endexample

In the next section, it will be shown that the conditions of
Theorem 2.2.3 hold when the underlying graph is of type $B_n$. Our
proof will also handle type $A_n$ as a special case.

\example{Example 2.2.5} Take the underlying graph $X$ to be of
type $D_n$ ($n\geq 4$). Let $\s_1,\s_2,\, \ldots,\s_n$ denote the
Coxeter generators, labelled so that $\s_3$ corresponds to the
branch node and $\s_1,\s_2$ commute with all generators except
$\s_3$. Consider the elements
$w=\s_2\s_3\s_4\s_3\s_1\s_2\s_3\notin W_c(D_n)$ and
$x=\s_1\s_2\s_4\s_3\in W_c(D_n)$.  Observe that $\s_1w>w$.  When
$C'_{\s_1}C'_w$ is written as a linear combination of
Kazhdan--Lusztig basis elements, the element $C'_x$ appears with
coefficient $1$. Thus, $x\leq_L w$, so that condition (vi) fails
when the underlying graph is of type $D_n$. Further computation
reveals that $x\sim_L w$, which shows that $W_c(D_n)$ is not a
union of left cells. This example also shows that condition (vi)
is violated (and that $W_c$ is not a union of left cells) in types
$E_6$, $E_7$ and $E_8$.  The incompatibility of $W_c$ with
Kazhdan--Lusztig cells in type $E$ was described explicitly in
\cite{{\bf 8}, \S9.9}.
\endexample

We remark that it is nevertheless true that the image under
$\theta$ of the set of all $C'_u$ indexed by $u\in W_c(D_n)$
equals the canonical basis of $\TL(D_n)$ (see \cite{{\bf 15},
Theorem 3.4}).  It is not known whether the corresponding
statement holds in type $E$.

\head 3. Type B \endhead

In this section, we study the compatibility of cells and fully
commutative elements when the underlying Coxeter graph $X$ is of
type $B_n$.

\subhead 3.1 Statement of results \endsubhead

Our main objective in \S3 is to prove the following

\proclaim{Theorem 3.1.1}  When the underlying Coxeter graph is of
type $B_n$, the equivalent conditions of Theorem 2.2.3 are
satisfied. In particular, the set $W_c(B_n)$ is closed under
$\geq_{LR}$ and so is a union of two-sided Kazhdan--Lusztig cells.
\endproclaim

This result is a strengthening of \cite{{\bf 12}, Theorem 2.2.1},
where it was shown that $\th(C'_w) \in \L$ for all $w \in W(B_n)$,
and thus that the basis in Lemma 2.2.1 for $\TL(B_n)$ agrees with
the canonical basis of $\TL(B_n)$.

Theorem 3.1.1 also implies the corresponding statement for type
$A_n$, which was previously known (see \cite{{\bf 3}, Proposition
3.1.1}).

\proclaim{Corollary 3.1.2} When the underlying Coxeter graph is of
type $A_n$, the equivalent conditions of Theorem 2.2.3 are
satisfied.
\endproclaim

\demo{Proof} Identify the Coxeter group $\,\,W(A_n)\,\,$ with the
parabolic subgroup of the Coxeter group $W(B_{n+1})$ that
corresponds to omission of the appropriate end generator.

Let $s \in S(A_n)$ and $w \in W(A_n)$ be such that $w \not\in
W_c(A_n)$ and $sw > w$.  By considering $C'_s C'_w$ and using
condition (vi) of Theorem 2.2.3 applied to type $B_{n+1}$, we see
that condition (vi) holds for type $A_n$. \qed\enddemo

The remaining cases arising from finite irreducible Coxeter groups
are $F_4$, $H_3$, and $H_4$. A series of computer calculations
using du Cloux's program ``Coxeter'' \cite{{\bf 1}} shows that
condition (vi) of Theorem 2.2.3 holds in each of these cases.

We can summarize our results as follows.

\proclaim{Corollary 3.1.3} Let $X$ be a Coxeter graph such that
$W(X)$ is finite and irreducible.  Then $W_c(X)$ is a union of
two-sided Kazhdan--Lusztig cells if and only if $X$ does not
contain $D_4$ as a subgraph. $\, \square$
\endproclaim

\subhead 3.2 Some combinatorial preparation \endsubhead

The following proposition describes a useful way to parse certain
reduced expressions. A proof can be found in \cite{{\bf 12}, Lemma
2.1.2}.

\proclaim{Proposition 3.2.1} Let $w\in W_c(B_n)$ and $s\in S(B_n)$
satisfy $ws\notin W_c(B_n)$.  There exists a unique $s'\in S(B_n)$
such that any reduced expression for $w$ can be parsed in one of
the following two ways.

\itemitem{\rm (i)}{$w=w_1sw_2s'w_3$, where $ss'$ has order
$3$, and $s$ commutes with every member of $c(w_2)\cup c(w_3);$}
\itemitem{\rm (ii)}{$w=w_1s'w_2sw_3s'w_4$, where
$ss'$ has order $4$, $s$ commutes with every member of $c(w_3)\cup
c(w_4)$, and $s'$ commutes with every member of $c(w_2)\cup
c(w_3)$.} $\, \square$
\endproclaim

The algebra $\TL(B_n)$ is known to be generated by the monomial
elements $b_s$, with $s$ ranging over all Coxeter generators,
subject to the following relations: $b_s^2=q_cb_s$, where $q_c =
[2]=v+v^{-1}$; $b_sb_{s'}=b_{s'}b_s$ if $s,s'$ commute;
$b_sb_{s'}b_s=b_s$ if $ss'$ has order 3;
$b_sb_{s'}b_sb_{s'}=2b_sb_{s'}$ if $ss'$ has order 4 (see
\cite{{\bf 9}, \S1}).

\proclaim{Lemma 3.2.2} Let $w\in W(B_n)$ and let $s_1s_2\cdots
s_m$ be a reduced expression for $w$.  Given integers $1\leq
i_1<i_2<\cdots <i_k\leq m$, we have $b_{s_{i_1}}b_{s_{i_2}}\cdots
b_{s_{i_k}}=aq_c^{\mu}b_{w'}$, where $a$ and $\mu$ are nonnegative
integers and $w'\in W_c(B_n)$.  Moreover, we have $w'\leq w$ and
$\ell(w')\leq k$.
\endproclaim

\demo{Proof} This follows by a simple induction on $k$, using the
subexpression characterization of Bruhat--Chevalley order together
with Proposition 3.2.1 and the relations for the monomial
generators given in the previous paragraph. \qed\enddemo

\remark{Remark 3.2.3}
We usually apply Lemma 3.2.2 in the following way.  Let $w=xyz$ be
a reduced product, and consider $b_x\te_yb_z$.  We would like to
know that this is a linear combination of monomial
basis elements $b_u$ with $u\leq w$.  To see that this is the case,
let $s_1s_2\cdots s_m$ be a reduced expression for $y$.
Then $b_x\te_yb_z$ equals
$$b_x\te_{s_1}\te_{s_2}\cdots
\te_{s_m}b_z=b_x(b_{s_1}-v^{-1})(b_{s_2}-v^{-1})\cdots
(b_{s_m}-v^{-1})b_z,$$ and we see that this expands into a
combination of $b_u$ with $u\leq w$ by Lemma 3.2.2.
\endremark

The next result gives useful information concerning the structure
constants for the monomial basis; it will be used repeatedly in
\S3.3.

\proclaim{Proposition 3.2.4} {\rm [{\bf 12}, Lemma 2.1.3] }Let $w\in
W_c(B_n)$ and let $s\in S(B_n)$.  We have
$b_wb_s=aq^{\mu}_cb_{w'}$ for some fully commutative $w'$ and some
nonnegative integers $a$ and $\mu$. Furthermore, one has {\rm (i)}
$\mu\leq 1;$ {\rm (ii)} $\ell(w's)<\ell(w');$ {\rm (iii)} $\mu=0$
if $\ell(ws')<\ell(w)$ for some $s'\in S(B_n)$ that does not
commute with $s$. $\, \square$
\endproclaim

Let $\s_1,\s_2,\ldots ,\s_n$ be the elements of $S(B_n)$, labelled
so that $\s_1\s_2$ has order 4 and $\s_i\s_{i+1}$ has order 3 for
all $i>1$.  Define, for each $1\leq r\leq n$, the set
$W^{(r)}=\{w\in W(B_r): i<r\Rightarrow \ell(\s_iw)>\ell(w)\}$.  It
is known that $W^{(r)}$ is a system of right coset representatives
for the parabolic subgroup $W(B_{r-1})$ of $W(B_r)$.  Moreover,
one has $\ell(xy)=\ell(x)+\ell(y)$ for all $x\in W(B_{r-1})$ and
$y\in W^{(r)}$ (see \cite{{\bf 13}, \S5.12}).  Thus, each $y\in W^{(r)}$
is the unique element of minimum length in the coset
$W(B_{r-1})y$. The elements of $W^{(r)}$ are given as follows:
$$\{ e,\, \s_r,\, \s_r\s_{r-1},\, \ldots ,\, \s_r\s_{r-1}\cdots
\s_2\s_1,\, \s_r\s_{r-1}\cdots \s_2\s_1\s_2,\, \ldots, $$
$$\s_r\s_{r-1}\cdots \s_2\s_1\s_2\cdots \s_{r-1}\s_r \}.$$ Note
that each element of $W^{(r)}$ has a unique reduced expression and
hence is fully commutative.

Any $w\in W(B_n)$ can be written uniquely as a product
$w_1w_2\cdots w_n$, where each $w_i\in W^{(i)}$.  By the previous
paragraph, this product is reduced.  Thus, if we delete each $w_i$
that equals the identity, and then replace each of the remaining
$w_i$ with its unique reduced expression, we obtain a ``normal"
reduced expression for $w$.

We frequently use without comment the following consequence of
the Exchange Condition, which is valid for any Coxeter system: if
$w\in W$ and $s\in S$, then $w$ has a reduced expression ending in
$s$ if and only if $\ell(ws)<\ell(w)$ (see \cite{{\bf 13}, \S5.8}).

\subhead 3.3 The $\A^-$-lattices $\L'_w$ \endsubhead

In the following series of lemmas, we study the $\A^-$-lattices in
$\TL(B_n)$ of the form $\L'_w$ (recall Definition 2.1.6).  Our
goal, which is accomplished in Proposition 3.3.10, is to prove
that $\te_w\in v^{-1}\L'_w$ whenever $w\notin W_c(B_n)$. This will
enable us to establish condition (v) of Theorem 2.2.3 for the case
where $X=B_n$.

\proclaim{Lemma 3.3.1}Let $x\in W_c(B_{n-1})$, let $w\in W^{(n)}$
and let $k\geq 0$. Suppose that there exist $s,s'\in S(B_n)$, with
$\ell(wss')<\ell(ws)<\ell(w)$, such that $b_x\te_{wss'}\in
v^{-k}\L'_{xwss'}$ and $b_x\te_{ws}\in v^{-k}\L'_{xws}$. Then
$b_x\te_w\in v^{-k}\L'_{xw}$.
\endproclaim

\demo{Proof} First note that when $b_x\te_w$ is written as a
linear combination of monomial basis elements $b_y$, the
coefficient of $b_y$ is nonzero only if $y\leq xw$ by Remark 3.2.3.
Thus, we may turn our attention to the degrees of the various coefficients.

By hypothesis, we can write $b_x\te_{wss'}$ as a sum of terms
$a_yb_y$ with each $a_y\in v^{-k}\A^-$.  Since $$b_x\te_{ws}=
b_x\te_{wss'}\te_{s'}=b_x\te_{wss'}(b_{s'}-v^{-1}),$$ we see that
$b_x\te_{ws}$ equals a sum of terms of the form
$a_yb_y(b_{s'}-v^{-1})$, where $a_y\in v^{-k}\A^-$.  We can use
this fact together with Proposition 3.2.4 (ii) and the hypothesis
on $b_x\te_{ws}$ to deduce that when $b_x\te_{ws}$ is written as a
linear combination of monomial basis elements $b_z$, either the
coefficient of $b_z$ lies in $v^{-k-1}\A^-$, or the coefficient
lies in $v^{-k}\A^-$ and $\ell(zs')<\ell(z)$.

Now consider the equalities
$$b_x\te_w=b_x\te_{ws}\te_s=b_x\te_{ws}(b_s-v^{-1}).$$ In view of
the previous paragraph, together with the fact that $s,s'$ do not
commute (since $w$ has a unique reduced expression), parts (i) and
(iii) of Proposition 3.2.4 enable us
to conclude that $b_x\te_w$ is a linear combination of monomial basis
elements $b_{z'}$ ($z'\leq xw$) with coefficients in $v^{-k}\A^-$, as
desired.
\qed\enddemo

\proclaim{Lemma 3.3.2}Let $x\in W_c(B_{n-1})$ and let $w\in
W^{(n)}$. Then $b_x\te_w\in \L'_{xw}$.
\endproclaim

\demo{Proof} We argue by induction on $\ell(w)$.  The lemma is
obviously true for $\ell(w)=0$, and if $\ell(w)=1$, then
$w=\s_n\notin c(x)$. Hence,
$b_x\te_w=b_x(b_{\s_n}-v^{-1})=b_{x\s_n}-v^{-1}b_x$, and one sees
that this last expression belongs to $\L'_{xw}$.

Suppose that $\ell(w)>1$.  There exist (uniquely determined)
Coxeter generators $s,s'$ such that $\ell(wss')<\ell(ws)<\ell(w)$.
By the inductive hypothesis, $b_x\te_{wss'}\in \L'_{xwss'}$ and
$b_x\te_{ws}\in \L'_{xws}$. But then $b_x\te_w\in \L'_{xw}$ by
Lemma 3.3.1 (taking $k=0$). The inductive step is complete.
\qed\enddemo

\proclaim{Proposition 3.3.3} We have $\te_w\in \L'_w$ for all
$w\in W(B_n)$.
\endproclaim

\demo{Proof} We proceed by induction on $\ell(w)$.  If
$\ell(w)=0$, then $w=e$ and we have $\te_e=b_e$. Suppose that
$\ell(w)>0$. Let $r>0$ be the smallest integer such that $w\in
W(B_r)$.  Write $w$ as a reduced product $w=yz$, where $y\in
W(B_{r-1})$ and $z\in W^{(r)}$.  We have $\te_w=\te_y\te_z$. By
the inductive hypothesis, we may write $\te_y$ as a linear
combination of monomial basis elements $b_x$ ($x\leq y$) with
coefficients in $\A^-$.  Thus, $\te_w$ equals a linear combination
of products of the form $b_x\te_z$ ($x\leq y$), with coefficients
in $\A^-$.

If we can show that any such product $b_x\te_z$ lies in $\L'_w$,
then the inductive step will be established.  But Lemma 3.3.2
gives us $b_x\te_z\in \L'_{xz}$, and the subexpression
characterization of Bruhat--Chevalley order gives $xz\leq w$. The
proof is complete. \qed\enddemo

The following two lemmas are needed to handle certain cases that
arise in the proofs of lemmas 3.3.6 and 3.3.8.

\proclaim{Lemma 3.3.4} Let $x\in W_c(B_{n-1})$ and let $w,w'\in
W^{(n)}$.  Suppose that $w=w'u$ (reduced) for some $u\in W(B_n)$.
Let $u'\in W(B_n)$ satisfy $u'\leq u$, and let $k\geq 1$. If
$\ell(u')<k$, then $v^{-k}b_x\te_{w'}b_{u'}\in v^{-1}\L'_{xw}$.
\endproclaim

\demo{Proof}   Fix a reduced expression $s_1s_2\cdots s_m$ for
$u'$.  We assume that $m<k$.  By Lemma 3.3.2, the product
$b_x\te_{w'}$ can be written as a sum of terms of the form
$a_yb_y$, where $a_y\in \A^-$ and $y\leq xw'$.  Thus,
$v^{-k}b_x\te_{w'}b_{u'}$ equals a sum of terms of the form
$$v^{-k}a_yb_yb_{u'}=v^{-k}a_yb_yb_{s_1}b_{s_2}\cdots b_{s_m},$$
where again $a_y\in \A^-$ and $y\leq xw'$.

By applying Proposition 3.2.4 (i) repeatedly ($m$ times on the
same term) and then applying Lemma 3.2.2, we find that
$v^{-k}a_yb_yb_{s_1}b_{s_2}\cdots b_{s_m}$ lies in
$v^{-1}\L'_{xw}$.\qed\enddemo

\proclaim{Lemma 3.3.5}Let $x\in W_c(B_{n-1})$ and let $w,w'\in
W^{(n)}$.  Suppose that $w=w'u$ (reduced) for some $u\in W(B_n)$.
Let $u'\in W(B_n)$ satisfy $u'\leq u$, and assume that $u'$ has a
unique reduced expression $s_1s_2\cdots s_m$. Then
$b_x\te_{w'}b_{u'}\in \L'_{xw}$ if there exists an $s\in S(B_n)$
that does not commute with $s_1$ and that satisfies either {\rm
(i)} $\ell(w's)<\ell(w');$ or {\rm (ii)} $\ell(xs)<\ell(x)$ and
$sw'=w's$.
\endproclaim

\demo{Proof} We first point out that by Remark 3.2.3, when
$b_x\te_{w'}b_{u'}$ is expressed as a linear combination of
monomial basis elements, all nonzero terms correspond to $y\in
W_c(B_n)$ that satisfy $y\leq xw$.

Suppose that (i) holds.  Write $w'=w''s$ (reduced).  We have
$b_x\te_{w'}\in \L'_{xw'}$ and $b_x\te_{w''}\in \L'_{xw''}$ by
Lemma 3.3.2. The equalities
$b_x\te_{w'}=b_x\te_{w''}\te_s=b_x\te_{w''}(b_s-v^{-1})$ together
with Proposition 3.2.4 (ii) imply that when $b_x\te_{w'}$ is
written as a sum of terms of the form $a_yb_y$ ($a_y\in \A^-$),
for each $y$, either $\ell(ys)<\ell(y)$ or else $a_y\in
v^{-1}\A^-$.

Thus, $b_x\te_{w'}b_{u'}$ is a sum of terms of the form
$a_yb_yb_{u'}=a_yb_yb_{s_1}b_{s_2}\cdots b_{s_m}$, with $y$ and
$a_y$ as described above.  If $\ell(ys)<\ell(y)$ then, since
$ss_1\neq s_1s$ and $s_is_{i+1}\neq s_{i+1}s_i$ for all $i<m$,
repeated applications of parts (ii) and (iii) of Proposition 3.2.4
give $a_yb_yb_{s_1}b_{s_2}\cdots b_{s_m}\in \L'_{xw}$.
On the other hand, if $a_y\in v^{-1}\A^-$, then
$a_yb_yb_{s_1}=a'b_{y'}$ with $a'\in \A^-$ by Proposition 3.2.4
(i).  Also, $\ell(y's_1)<\ell(y')$ by part (ii) of the same
proposition. Since $s_is_{i+1}\neq s_{i+1}s_i$ for all $i<m$, we
have $a_yb_yb_{s_1}b_{s_2}\cdots b_{s_m}=a'b_{y'}b_{s_2}\cdots
b_{s_m}\in \L'_{xw}$ (again by repeated
applications of parts (ii) and (iii) of Proposition 3.2.4).

Suppose now that (ii) holds.  Write $x=x's$ (reduced).  Since
$sw'=w's$, we have $b_x\te_{w'}=b_{x'}\te_{w'}b_s$; hence, by
Proposition 3.2.4 (ii), when $b_x\te_{w'}$ is written as a sum of
terms of the form $a_yb_y$, we have $\ell(ys)<\ell(y)$ whenever
$a_y\neq 0$.  But then the reasoning from the previous paragraph
gives $a_yb_yb_{s_1}b_{s_2}\cdots b_{s_m}\in \L'_{xw}$.
\qed\enddemo

The following lemma and Lemma 3.3.8 are needed to handle the
inductive step of Lemma 3.3.9.

\proclaim{Lemma 3.3.6} Let $x\in W_c(B_{n-1})$ and let $w\in
W^{(n)}$. Suppose that there exists $s\in S(B_n)$, with
$\ell(ws)<\ell(w)$, such that $xws$ is fully commutative but $xw$
is not. Then $b_x\te_w\in v^{-1}\L'_{xw}$.
\endproclaim

\demo{Proof} By the hypothesis of the lemma, together with
Proposition 3.2.1, there are two possibilities concerning the
nature of $w$ and $x$: Either (a) $w=\s_n\s_{n-1}\cdots
\s_{r+1}\s_r$ ($2 \leq r < n$) and $\ell(x\s_r)<\ell(x)$; or (b)
$w=\s_n\s_{n-1}\cdots \s_2\s_1\s_2$ and $\ell(x\s_1)<\ell(x)$.

Suppose that case (a) holds.  Let $w'=w\s_r\s_{r+1}$.  We have $$
b_x\te_w=b_{x\s_r}b_{\s_r}\te_{w'}\te_{\s_{r+1}}\te_{\s_r} =
b_{x\s_r}\te_{w'}b_{\s_r}(b_{\s_{r+1}}-v^{-1})
(b_{\s_r}-v^{-1}).$$ This last expression expands into the
product $$ b_{x\s_r}\te_{w'}(b_{\s_r}b_{\s_{r+1}}b_{\s_r} -
v^{-1}b_{\s_r}b_{\s_{r+1}}-v^{-1}b_{\s_r}^2+v^{-2}b_{\s_r}).$$
Using the relations for the monomial generators given in
\S3.2, we may simplify this last expression to $$
b_{x\s_r}\te_{w'}(-v^{-1}b_{\s_r}b_{\s_{r+1}})=
b_x\te_{w'}(-v^{-1}b_{\s_{r+1}}).$$ Thus,
$b_x\te_w=-v^{-1}b_x\te_{w'}b_{\s_{r+1}}$.  This last
expression is easily seen to lie in $v^{-1}\L'_{xw}$ by Lemma
3.3.5 (ii) (taking $s = \s_r$).

We turn to case (b), wherein $w=\s_n\s_{n-1}\cdots \s_2\s_1\s_2$
and $x$ has a reduced expression ending in $\s_1$. Let
$w''=w\s_2\s_1\s_2$.  Observe that $b_x\te_w$ equals
$$b_{x\s_1}b_{\s_1}\te_{w''}\te_{\s_2}\te_{\s_1}\te_{\s_2} =
b_{x\s_1}\te_{w''}b_{\s_1}(b_{\s_2}-v^{-1})(b_{\s_1}-v^{-1})
(b_{\s_2}-v^{-1}).$$ By expanding and then simplifying this last
expression, we obtain
$$b_{x\s_1}\te_{w''}(-v^{-1}b_{\s_1}b_{\s_2}b_{\s_1}+
v^{-1}b_{\s_1})=b_x\te_{w''}(-v^{-1}b_{\s_2\s_1}+v^{-1}).$$ Thus,
$b_x\te_w=-v^{-1}b_x\te_{w''}b_{\s_2\s_1}+v^{-1}b_x\te_{w''}$.
This last pair of terms lies in $v^{-1}\L'_{xw}$ by
Lemma 3.3.5 (ii).
\qed\enddemo

The next lemma is required to address a certain term that arises
in the proof of Lemma 3.3.8.

\proclaim{Lemma 3.3.7} Let $x\in W_c(B_{n-1})$ and let
$w=\s_n\s_{n-1}\cdots \s_r\s_{r-1}$, with $2\leq r<n$. Let
$w'=w\s_{r-1}\s_r\s_{r+1}$. Suppose that $\ell(x\s_r)<\ell(x)$.
Then $b_x\te_{w'}b_{\s_{r-1}\s_{r+1}}\in \L'_{xw}$.
\endproclaim

\demo{Proof} Note first that when
$b_x\te_{w'}b_{\s_{r-1}\s_{r+1}}$ is written as a linear
combination of monomial basis elements $b_y$, every nonzero term
is indexed by a $y$ satisfying $y\leq xw$ by Remark 3.2.3.
If $w'=e$, then $b_x\te_{w'}b_{\s_{r-1}\s_{r+1}}=
b_xb_{\s_{r-1}}b_{\s_{r+1}}$, which equals $ab_yb_{\s_{r+1}}$
for some integer $a$ and some $y\in W_c(B_{n-1})$ by
Proposition 3.2.4 (iii). Since $w'=e$, we have
$\s_{r+1}=\s_n\notin c(y)$, hence the product $y\s_{r+1}$ is
reduced and fully commutative.  Therefore,
$ab_yb_{\s_{r+1}}=ab_{y\s_{r+1}}\in \L'_{xw}$.

Now assume $w'\neq e$.  Let $w''=w'\s_{r+2}$ and let $x'=x\s_r$.
By Lemma 3.3.2, we can write $b_x\te_{w'}$ as a sum of terms of
the form $a_yb_y$, where $y\leq xw'$ and $a_y\in \A^-$.  Moreover,
we have $\ell(y\s_r)<\ell(y)$ whenever $a_y\neq 0$, owing to the
equality $b_x\te_{w'}=b_{x'}\te_{w'}b_{\s_r}$ and Proposition
3.2.4 (ii).  Thus,
$b_x\te_{w'}b_{\s_{r-1}\s_{r+1}}$ is a sum of terms of the form
$a_yb_yb_{\s_{r-1}\s_{r+1}}$, where $y\leq xw'$ and
$\ell(y\s_r)<\ell(y)$ whenever $a_y\neq 0$.

We can say more about the terms $a_yb_yb_{\s_{r-1}\s_{r+1}}$.
Since $b_x\te_{w'}=b_x\te_{w''}\te_{\s_{r+2}}=b_x\te_{w''}
(b_{\s_{r+2}}-v^{-1})$, we see (as in the proof of Lemma 3.3.1)
that for each $y\leq xw'$, either
$a_y\in v^{-1}\A^-$ or else $a_y\in \A^-$ and
$\ell(y\s_{r+2})<\ell(y)$. In the case where $a_y\in v^{-1}\A^-$
and $a_y\neq 0$, we have $a_yb_yb_{\s_{r-1}}=a'b_{y'}$ with
$a' \in v^{-1}\A^-$ by Proposition 3.2.4 (iii) (here, we are
using the fact that $\ell(y\s_r)<\ell(y)$). But then
$a_yb_yb_{\s_{r-1}}b_{\s_{r+1}}=
a'b_{y'}b_{\s_{r+1}}=a''b_{y''}$, with $a''\in
\A^-$ by Proposition 3.2.4 (i).

Consider now the other case, where $a_y\in \A^-$, $a_y\neq 0$ and
$\ell(y\s_{r+2})<\ell(y)$. Here, we may write $y$ as a fully
commutative reduced product $y=y_1\s_r\s_{r+2}$. Notice that
$\ell(y\s_{r-1})>\ell(y)$.  If $y\s_{r-1}$ is fully commutative,
then, since $y\s_{r-1}=y_1\s_r\s_{r-1}\s_{r+2}$ and
$\s_{r+1}\s_{r+2}\neq \s_{r+2}\s_{r+1}$, we have
$a_yb_yb_{\s_{r-1}}b_{\s_{r+1}}=a_yb_{y\s_{r-1}}b_{\s_{r+1}}=
a'b_{y'}$ for some $a'\in \A^-$ by Proposition 3.2.4 (iii).

If $y\s_{r-1}$ is not fully commutative, then there are two
possibilities to consider, depending on whether $r-1$ is greater
than or equal to $1$. By Proposition 3.2.1 we have, after applying
commutations to obtain a suitable reduced expression for $y$ if
necessary, either (a) $y=y_2\s_{r-1}\s_r\s_{r+2}$ (reduced) if
$r-1>1$; or (b) $y=y_3\s_r\s_{r-1}\s_r\s_{r+2}$ (reduced) if
$r-1=1$. The argument for (b) is essentially the same as that for
(a), so we treat only (a).

The product $a_yb_yb_{\s_{r-1}}b_{\s_{r+1}}$ equals
$a_yb_{y_2}b_{\s_{r-1}}b_{\s_r}b_{\s_{r+2}}b_{\s_{r-1}}
b_{\s_{r+1}}$. Since $\s_{r+2}$ and $\s_{r-1}$ commute, this last
expression equals $a_yb_{y_2}b_{\s_{r-1}}b_{\s_r}b_{\s_{r-1}}
b_{\s_{r+2}}b_{\s_{r+1}}=a_yb_{y_2}b_{\s_{r-1}}b_{\s_{r+2}}
b_{\s_{r+1}}$. Now, recall that $y=y_2\s_{r-1}\s_r\s_{r+2}$ is a
fully commutative reduced product.  Since $\s_{r+2}$ and $\s_r$
commute, the product $y_2\s_{r-1}\s_{r+2}$ is reduced and fully
commutative.  It follows that
$a_yb_yb_{\s_{r-1}}b_{\s_{r+1}}=a_yb_{y_2}b_{\s_{r-1}}
b_{\s_{r+2}}b_{\s_{r+1}}$ equals
$a_yb_{y_2\s_{r-1}\s_{r+2}}b_{\s_{r+1}}$, and this last expression
equals $a''b_{y''}$ for some $a''\in \A^-$ by Proposition 3.2.4
(iii). \qed\enddemo

\proclaim{Lemma 3.3.8} Let $x\in W_c(B_{n-1})$ and let $w\in
W^{(n)}$. Suppose that there exist $s,s'\in S(B_n)$, with
$\ell(wss')<\ell(ws)<\ell(w)$, such that $xwss'$ is fully
commutative but $xws$ is not. Then $b_x\te_w\in v^{-1}\L'_{xw}$.
\endproclaim

\demo{Proof} As in the proof of Lemma 3.3.6, we can use
Proposition 3.2.1 to divide the argument into two cases:
Either (a)
$w=\s_n\s_{n-1}\cdots \s_{r+1}\s_r\s_{r-1}$, with $2 \leq r < n$ and
$\ell(x\s_r)<\ell(x)$; or (b) $n \geq 3$ and $w=\s_n\s_{n-1}\cdots
\s_2\s_1\s_2\s_3$ with $\ell(x\s_1)<\ell(x)$.

Suppose first that (a) holds.  Let $w'=w\s_{r-1}\s_r\s_{r+1}$ and
let $x'=x\s_r$. The product $b_x\te_w$ equals
$b_{x'}b_{\s_r}\te_{w'}\te_{\s_{r+1}}\te_{\s_r} \te_{\s_{r-1}}$,
which in turn can be written as
$$b_{x'}\te_{w'}b_{\s_r}(b_{\s_{r+1}}-v^{-1})
(b_{\s_r}-v^{-1})(b_{\s_{r-1}}-v^{-1}).$$ After expanding and then
simplifying this last expression, we obtain
$$b_{x'}\te_{w'}(-v^{-1}b_{\s_r}b_{\s_{r+1}}
b_{\s_{r-1}}+v^{-2}b_{\s_r}b_{\s_{r+1}})=b_x\te_{w'}
(-v^{-1}b_{\s_{r+1}\s_{r-1}}+v^{-2}b_{\s_{r+1}}).$$

Thus, $b_x\te_w=-v^{-1}b_x\te_{w'}b_{\s_{r+1}\s_{r-1}}+
v^{-2}b_x\te_{w'}b_{\s_{r+1}}$.  The first of these two terms lies
in $v^{-1} \L'_{xw}$ by Lemma 3.3.7 and the second lies in $v^{-1}
\L'_{xw}$ by Lemma 3.3.4.

Suppose now that (b) holds.  Let $w'=w\s_3\s_2\s_1\s_2\s_3$ and
let $x'=x\s_1$.  The product $b_x\te_w$ equals
$b_{x'}b_{\s_1}\te_{w'}\te_{\s_3}\te_{\s_2}\te_{\s_1}\te_{\s_2}
\te_{\s_3},$ which in turn equals
$$b_{x'}\te_{w'}b_{\s_1}(b_{\s_3}-v^{-1})(b_{\s_2}-v^{-1})
(b_{\s_1}-v^{-1})(b_{\s_2}-v^{-1})(b_{\s_3}-v^{-1}).$$  This last
expression expands and then simplifies to
$$b_x\te_{w'}(v^{-2}b_{\s_3\s_2\s_1}+v^{-2}b_{\s_2\s_1\s_3}-
2v^{-2}b_{\s_3}-v^{-3}b_{\s_2\s_1}+v^{-3}).$$  The products
$-2v^{-2}b_x\te_{w'}b_{\s_3}$, $-v^{-3}b_x\te_{w'}b_{\s_2\s_1}$,
and $v^{-3}b_x\te_{w'}$ belong to $v^{-1} \L'_{xw}$ by Lemma
3.3.4. Consider the product $v^{-2}b_x\te_{w'}b_{\s_3\s_2\s_1}$.
If $w'=e$, then $n=3$ and $x\s_3\s_2\s_1$ is fully commutative;
the latter holds because (1) $x$ is fully commutative and has no
reduced expression ending in $\s_2$ (as $\ell(x\s_1)<\ell(x)$ and
$\s_1$, $\s_2$ do not commute) and (2) the presence of $\s_3$
precludes the possibility of obtaining a substring of the form
$\s_1\s_2\s_1\s_2$ or $\s_2\s_1\s_2\s_1$ through commutation
moves. Hence,
$v^{-2}b_x\te_{w'}b_{\s_3\s_2\s_1}=v^{-2}b_{x\s_3\s_2\s_1}\in
v^{-1}\L'_{xw}$.  If $w'\neq e$, then
$v^{-2}b_x\te_{w'}b_{\s_3\s_2\s_1}\in \L'_{xw}$ by Lemma 3.3.5 (i)
(taking $s=\s_4$).  Finally, we have $b_x \te_{w'} b_{\s_2\s_1}
\in \L'_{xw'\s_3\s_2\s_1}$ by Lemma 3.3.5 (ii) (taking $s =
\s_1$), which, by Proposition 3.2.4 (i), shows that
$v^{-2}b_x\te_{w'}b_{\s_2\s_1\s_3}\in v^{-1} \L'_{xw}$.
\qed\enddemo

\proclaim{Lemma 3.3.9} Let $x\in W_c(B_{n-1})$ and let $w\in
W^{(n)}$. Suppose that $xw$ is not fully commutative. Then
$b_x\te_w\in v^{-1}\L'_{xw}$.
\endproclaim

\demo{Proof} We proceed by induction on $\ell(w)$.  The lemma is
(vacuously) true for $\ell(w)<2$, since $xw$ is fully commutative
for such $w\in W^{(n)}$.  Suppose that $\ell(w)\geq 2$. Let $s,s'$
be (the) Coxeter generators that satisfy
$\ell(wss')<\ell(ws)<\ell(w)$.  If $xwss'$ is not fully
commutative, then neither is $xws$; by the inductive hypothesis,
$b_x\te_{wss'}\in v^{-1}\L'_{xwss'}$ and $b_x\te_{ws}\in
v^{-1}\L'_{xws}$. But then $b_x\te_w\in v^{-1}\L'_{xw}$ by Lemma
3.3.1 (taking $k=1$).

On the other hand, suppose that $xwss'$ is fully commutative.  If
$xws$ is not fully commutative, then Lemma 3.3.8 gives us
$b_x\te_w\in v^{-1}\L'_{xw}$. Finally, if $xws$ is fully
commutative, then Lemma 3.3.6 applies.

The inductive step is complete. \qed\enddemo

\proclaim{Proposition 3.3.10}
We have $\te_w\in v^{-1}\L'_w$ for all
$w\in W(B_n)\backslash W_c(B_n)$.
\endproclaim

\demo{Proof}  Let $w\in W(B_n)\backslash W_c(B_n)$.
Write $w=w_1w_2\cdots w_n$, where each $w_i\in W^{(i)}$.  There
exists a unique integer $r>1$ such that $w_1w_2\cdots w_{r-1}$ is
fully commutative and $w_1w_2\cdots w_r$ is not.  Let
$y=w_1w_2\cdots w_{r-1}$.  Given the nature of the
representative $w_r$, it follows from Proposition 3.2.1 that
$y=y's$ (reduced), where $s\in S(B_{r-1})$ is such that $sw_r$ is
not fully commutative. By Proposition 3.3.3, we have $\te_{y'}\in
\L'_{y'}$ and $\te_y\in \L'_y$. Combining this with the equality
$\te_y=\te_{y'}(b_s-v^{-1})$, we find that if $\te_y$ is written
as a linear combination of monomial basis elements $b_x$ ($x\leq
y$), then for each $x$, either the coefficient of $b_x$ lies in
$v^{-1}\A^-$, or the coefficient lies in $\A^-$ and
$\ell(xs)<\ell(x)$.

Thus, $\te_{yw_r}=\te_y\te_{w_r}$ equals a sum of terms of the
form $a_xb_x\te_{w_r}$ ($x\leq y$), where for each
$x$, either $a_x\in v^{-1}\A^-$, or else $a_x\in \A^-$ and $xw_r$
is not fully commutative. Now, by Lemma 3.3.2, we have
$b_x\te_{w_r}\in \L'_{xw_r}$. Therefore, in the case where $a_x\in
v^{-1}\A^-$, we have $a_xb_x\te_{w_r}\in v^{-1}\L'_{xw_r}$. In
the case where $a_x\in \A^-$ and $xw_r\notin W_c$, Lemma 3.3.9
gives us $a_xb_x\te_{w_r}\in v^{-1}\L'_{xw_r}$. We have so far
shown that $\te_{yw_r}\in v^{-1}\L'_{yw_r}$.

The product $\te_{yw_r}\te_{w_{r+1}}$ must then lie in
$v^{-1}\L'_{yw_rw_{r+1}}$. To see this, expand $\te_{yw_r}$ into a
sum of terms $a'_zb_z$ ($z\leq yw_r$), where each $a'_z \in
v^{-1}\A^-$. Then $\te_{yw_r}\te_{w_{r+1}}$ equals a sum of terms
$a'_zb_z\te_{w_{r+1}}$, each of which must lie in
$v^{-1}\L'_{yw_rw_{r+1}}$ by Lemma 3.3.2 (here, $r+1$ is playing
the role of $n$ in the lemma). By iterating this argument, we are
able to conclude that $\te_{w}=\te_{yw_r}\te_{w_{r+1}}\cdots
\te_{w_n}\in v^{-1}\L'_w$, as desired.\qed\enddemo

\subhead 3.4 Proof of Theorem 3.1.1 \endsubhead

We are now in a position to give a proof that the conditions of
Theorem 2.2.3 are satisfied when the underlying Coxeter graph is
of type $B_n$.

\demo{Proof of Theorem 3.1.1} We shall verify condition (v) of
Theorem 2.2.3.  It is required to prove that for each $w\in
W(B_n)\backslash W_c(B_n)$, the element $\theta(\T_w)=\te_w$ lies
in $v^{-1}\L$.  Proposition 3.3.10 gives $\te_w\in v^{-1}\L'_w$
for such $w$.  Thus, the theorem will follow if we can show that
$\L'_w=\L_w$.

Fix an arbitrary $w$.  Proposition 3.3.3 gives $\L'_w\supseteq
\L_w$. Let $x \in W_c$ satisfy $x\leq w$. By Proposition 3.3.3, we
may write $\te_x = \sum a_y b_y$, where for all $y$ we have $y
\leq x$, $y \in W_c$ and $a_y \in \A^-$. Moreover, given any
reduced expression $s_1s_2\cdots s_m$ for $x$, we have
$$\te_x=\te_{s_1}\te_{s_2}\cdots \te_{s_m}=
(b_{s_1}-v^{-1})(b_{s_2}-v^{-1})\cdots (b_{s_m}-v^{-1}),$$ from
which we can see that $a_x = 1$. It now follows by a
straightforward induction on the Bruhat--Chevalley order that any
$b_x$ with $x\leq w$ can be written as an $\A^-$-linear
combination of basis elements $\te_y$ with $y\leq x$.  That is, we
have $\L'_w\subseteq \L_w$, as well. \qed\enddemo

\remark{Remark 3.4.1} Another possible approach to proving Theorem
3.1.1 would be to use the combinatorial classification of cells
achieved by Garfinkle in \cite{{\bf 5}, {\bf 6}, {\bf 7}}.
\endremark

We conclude with a discussion of an application: it is possible as
a consequence of Theorem 3.1.1 to describe very explicitly the
structure of the fully commutative left, right and two-sided cells
in type $B$.  By \cite{{\bf 12}, Theorem 2.2.1}, the canonical basis
for $\TL(B_n)$ is the image under $\th$ of the set $\{C'_w: w \in
W_c(B_n)\}$, and by condition (iii) of Theorem 2.2.3, all other
$C'_w$ are mapped to zero.  The diagram calculus for the canonical
basis of $\TL(B_n)$ given in \cite{{\bf 10}, Theorem 2.2.5} can now be
used to describe the various kinds of fully commutative
Kazhdan--Lusztig cells.

By \cite{{\bf 16}, Theorem 1.10}, each
left (respectively, right) cell in a crystallographic Coxeter group
$W$, such as $W(B_n)$, contains a unique distinguished involution
$d$.  In type $A$, any left cell and right cell from the same
two-sided cell intersect in a single element.  Our approach can be
used to describe the situation in type $B$ in an elementary way.

Let $n \geq 2$, and let $W' \cong W(A_{n-1})$ be the parabolic
subgroup of $W(B_n)$ obtained by omitting the first generator. Let
$d, d' \in W(B_n)$ be distinguished involutions in the same fully
commutative two-sided cell, $I_T$.  Let $I_R$ be the right cell
containing $d$, and let $I_L$ be the left cell containing $d'$.
Let $k = |I_R \cap I_L|$.  Then we have $k = 1$ if exactly one of
the elements $d, d'$ lies in $W'$, or if $I_T \subseteq W'$, or if
$I_T \cap W' = \emptyset$.  Otherwise, we have $k = 2$.

\leftheadtext{} \rightheadtext{}

\end